\magnification=1200 %
\parindent 25pt %
\parskip=0pt %
\baselineskip=17pt %

\input amssym.def
\input amssym

\vbox to .52in { }
\centerline{\bf On the semi-simplicity of the quantum cohomology  
}
\centerline{\bf algebras of complete intersections}
\par
\vskip .2in
\centerline{Gang Tian and Geng Xu}
\par 
\vskip .4in
\centerline{\bf 0. Introduction}
\par
\vskip .1in
{In this short note, we study the semi-simplicity of the quantum cohomology 
algebras of smooth Fano complete intersections in projective space.}
\par
{Let $V$ be a smooth Fano manifold with a symplectic form $\omega $. 
The quantum cohomology on $V$ is the cohomology $H^*(V,\Bbb Z \{H_2(V)\})$ 
with a ring structure defined by the GW-invariants. The Novikov ring 
$\Bbb Z\{H_2(V)\}$ 
 can be 
described as follows: choose a basis
$q_1,\cdots ,q_s$ of $H_2(V,\Bbb Z)$, we identify the monomial
$q^d=q_1^{d_1}\cdots q_s^{d_s}$ with the sum $\sum _{i=1}^s d_i q_i$.
This turns $H_2(V)$ into a multiplicative ring,
i.e., $q^d \cdot q^{d'} = q^{d+d'}$. 
This multiplicative ring has a natural grading defined by
$\deg (q^d) = 2 c_1(V) (\sum d_i q_i)$. Then 
$\Bbb Z\{H_2(V)\}$ is the graded homogeneous ring
generated by all formal power series
$\sum_{d=(d_1,\cdots,d_s)} n_d q^d$ satisfying: $n_d \in \Bbb Z$, 
all $q^d$ with $n_d \not= 0$ have the same degree,    
and the number of $n_d$ with $ \omega (\sum d_i q_i) \le c$ is finite
for any $c>0$.} 
\par
{Now we can define a ring structure on $H^*(V,\Bbb Z\{ H_2(V)\} )$.
For any $\alpha ^*$, $\beta ^*$ in $H^*(V,\Bbb Z)$, we define
the quantum multiplication $\alpha ^*\bullet \beta ^*$ by}  
$$\alpha ^*\bullet \beta ^* (\gamma ) = 
\sum _{A \in H_2(V,\Bbb Z)} \Psi^V_{(A,0,3)}
(\alpha, \beta , \gamma)q^A,$$  
{where $\gamma \in H_*(V,\Bbb Z)$, and  $\Psi^V_{(A,0,3)}(\alpha, \beta , \gamma)$ 
are the Gromov-Witten invariants (cf. [RT]). } 
\par
{In fact, there is a family of quantum multiplications. 
Let $\{\beta _a\}_{1\le a\le L}$ be an integral basis of 
$H_*(V,\Bbb Z)$ modulo torsions.
Any $w \in  H^*(V,\Bbb C)$ can be written as
$\sum t_a \beta^* _a$. Clearly, $w\in H^*(V,\Bbb Z)$ if all $t_a$ are integers.
We define the quantum multiplication $\bullet _w$ by} 
$$
\alpha ^* \bullet _w \beta ^* (\gamma )
 =  \sum _A \sum _{k\ge 0} {\epsilon (\{a_i\}) \over k!}
\Psi^V_{(A,0,k+3)}(\alpha, \beta , \gamma, \beta _{a_1}, \cdots, \beta _{a_k})
t_{a_1}\cdots t_{a_k} \, q^A $$
{where $\alpha , \beta , \gamma \in H_*(V,\Bbb Z)$,
 $\epsilon (\{a_i\})$ is the sign of the induced permutation
on odd dimensional $\beta _a$, and 
$\Psi^V_{(A,0,k+3)}(\alpha, \beta , \gamma, 
\beta _{a_1}, \cdots, \beta _{a_k})$ 
are the Gromov-Witten invariants.
Obviously, this multiplication reduces to $\bullet$ at $w=0$.
It was shown in [RT] that the quantum multiplications $\bullet _w$ are 
associative. }
\par 
\par 
{Now we restrict ourselves to the case $w\in W=H^{\hbox{even}}
(V,\Bbb C)$. In this case, 
 Dubrovin [D] observed that the quantum multiplications 
$\bullet_w$ induce the structure of a  Frobenius algebra on $W$.
 For any $w\in W, w=\sum t_a\beta^*_a$, put }
$$X(w) = \sum _{a\leq N} (\deg(\beta _a^*)-2) t_a \partial _a w - 2 c_1(V) 
;$$
{here we arrange the basis $\{ \beta_a\}$ so that $\beta_a$ is an even class if 
and only if $a\leq N$. 
We say that $V$ is semi-simple in the sense of Dubrovin, if for a 
generic $w$, the quantum multiplication $X(w)\bullet_w$ on $H^*(V, \Bbb C)$ 
has only simple eigenvalues. It is conjectured (cf. [T]) that
any Fano manifold is semi-simple in the 
sense of Dubrovin. The question we are interested in is the 
following conjecture [T]:}
\par
\vskip .1in
\noindent
{CONJECTURE (Tian).} {\it Any Fano manifold is semi-simple in the
sense of Dubrovin.}
\par
\vskip .1in
{There is a weaker version of the above conjecture. Let
$H^*_{inv}(V, \Bbb C)$ be the subring of $H^*(V, \Bbb C)$ with 
the cup product, generated by $H^2(V, \Bbb C)$. For any algebraic
variety $V$, it is believed
(cf. [T]) that for any $ w \in H^*_{inv}(V, \Bbb C)$,
the quantum multiplication $\bullet_w$
preserves the subspace $H^*_{inv}(V, \Bbb C)$, namely, if
$\alpha, \beta $ are in $ H^*_{inv}(V, \Bbb C)$, 
so is $\alpha \bullet_w\beta$.
Let us define $X_{inv}(w)$ be the restriction of $X$ to 
$H^*_{inv}(V, \Bbb C)$. If the above belief is true,
then $X_{inv}$ acts on $H^*_{inv}(V, \Bbb C)$. It is conjectured
in [T] that if $V$ is Fano, then $X(w)\bullet_w$ is semi-simple 
for a generic $w$ in $H^*_{inv}(V, \Bbb C)$.}
\par
{Now, let $V\subset {\bf P^{n+r}}$ be a Fano complete intersection of 
dimension $n\geq 3$. 
Then the ordinary cohomology algebra $H^*(V,\Bbb C)$ is 
generated by the hyperplane class $H$ and the primitive cohomology 
$H^n(V,\Bbb C)_o$. 
In particular, $H^*_{inv}(V, \Bbb C)$ be the subspace 
of $H^*(V, \Bbb C)$ 
generated by the hyperplane class $H$. It was observed in
[T] that for any $ w \in H^*_{inv}(V, \Bbb C)$,
the quantum multiplication $\bullet_w$
preserves the subspace $H^*_{inv}(V, \Bbb C)$.
Then $X_{inv}$ acts on $H^*_{inv}(V, \Bbb C)$ for any $w$
in $H^*_{inv}(V, \Bbb C)$.
An important case of the above conjecture is
that $X(w)\bullet_w$ is semi-simple for a generic $w$ in 
$H^*_{inv}(V, \Bbb C)$ and any complete intersection $V$ (cf. [T]).  
the result  we have  is the following: }
\par
\vskip .1in
\noindent
{THEOREM 1.}  
{\it Let $V \subset {\bf P^{n+r}} (n\geq 3)$ be a smooth complete 
intersection  of degree $(d_1,d_2, \cdots, d_r)$. 
If $n > 2\sum_{i=1}^r (d_i-1)-1$ (except $n=7$ and $\sum_{i=1}^r (d_i-1) 
= 2$), 
then $X(w) \bullet_{w}$, 
which acts on $H^*_{inv}(V, \Bbb C)$, 
is semi-simple in the sense of Dubrovin for a generic 
$w\in H^*_{inv}(V, \Bbb C)$.} 
\par
\vskip .1in
{Here we use Beauville's computation of quantum cohomology algebra
for complete intersections satisfying the condition in the above 
theorem (cf. [B]).
It seems that this condition on degree can be removed by using 
recent results in
[G].}
\par
{Our motivation for the study of semi-simplicity comes from the following 
result of Dubrovin [D]: if $\{ \alpha_1, \cdots, \alpha_m \} $ is a basis 
of $H^*(V,\Bbb C)$ such that the quantum multiplication $X(w)\bullet_w $ on 
$H^*(V,\Bbb C)$ with respect to this basis has only simple eigenvalues, then 
the integrable system defined by
using the Gromov-Witten prepotential $\Phi^V$ 
(see section 2) can be extented meromorphically to $({\bf P^1})^m$.}
\par 
{The main tool of this paper is an elementary lemma in section 1, the 
rest are standard computations. The method should also apply to 
some other Fano complete intersections. In this note, by a rational 
curve on $V$ we mean a simple genus 0 J-holomorphic curve for some 
generic almost complex structure $J$ on $V$.}
\par  
{Throughout this paper we work over the complex number field $\Bbb C$.} 

\par
\vskip .2in
\centerline{\bf 1. An algebraic lemma}
\par
\vskip .1in
{The main tool of our computations is the following elementary lemma.}
\par
\vskip .1in
\noindent
{LEMMA 1.} {\it Assume that }
$$g(y,z) = y^m + g_1(z) y^{m-1} + \cdots + g_m(z)$$
{\it is a polynomial of $y$, and $g_1(z), \cdots, g_m(z)$ are holomorphic 
functions of $z=(z_1, z_2, \cdots, z_N)$ defined in a neighborhood of 
$0=(0, \cdots, 0)\in \Bbb C^N$. If the only repeated root of the 
polynomial $g(y,0)=0$ is $y=0$, and the polynomial $g(y,z)=0$ does not 
have distinct roots for generic $z\in \Bbb C^N$, then}
$$g_m(z) = 0 \quad \hbox{ mod } (z_1, z_2, \cdots, z_N)^2,$$
{\it that is, the constant term and linear terms of $g_m(z)$ are all 0.}
\par
\vskip .1in
\noindent
{\it Proof.} 
{Let ${\cal O}_0$ denote the ring of holomorphic functions defined in 
some neighborhood of $(0,0)\in \Bbb C^N \times \Bbb C$. Since all the 
repeated roots of $g(y,z)=0$ appear nearby $y=0$, we will study 
$g(y,z)$ in ${\cal O}_0$. By ([GH] page 8-11), ${\cal O}_0$ is a unique 
factorization domain module units, we can write}
$$g(y,z) = (f_1(y,z))^{i_1} \cdots (f_h(y,z))^{i_h},$$
{where $f_i(y,z)\in {\cal O}_0$ are 
distinct and irreducible Weierstrass polynomials in $y$.}  
\par
{Now we claim that $i_j > 1$ for some $j$. Otherwise, let $V_j=\{ f_j(y,z)=0\}$,
 then }
$$ g(y,z)=f_1(y,z) \cdots f_h(y,z), \quad \{ g(y,z)=0 \} = V_1 \cup \cdots \cup V_h.$$
{Since the Weierstrass polynomial $f_i(y,z)$ is irreducible, and 
${{\partial f_i}\over {\partial y}} (y,z)$ has degree lower than $f_i(y,z)$ in $y$, 
we conclude that 
 $f_i(y,z)$ and ${{\partial f_i}\over {\partial y}} (y,z)$ are relative prime 
in ${\cal O}_0$. Hence the analytic variety }
$$\{ f_i(y,z) = 0\} \cap \bigl\{ {{\partial f_i}\over {\partial y}} (y,z) = 0\bigr\} 
\subset \Bbb C^N \times \Bbb C$$
{has dimension $N-1 < \hbox{dim }\Bbb C^N$. 
Therefore the polynomial $f_i(y,z) = 0$ has 
distinct roots for generic $z\in \Bbb C^N$.}
\par
{Similarly, the analytic variety $\{ f_i = 0\} \cap \{ f_j = 0\} 
\subset \Bbb C^N \times \Bbb C$ has 
dimension $N-1<N$ for $i\not= j$ by Weak Nullstellensatz ([GH] p. 11). 
Hence the polynomials} 
$$f_i(y,z)=0, \quad f_j(y,z)=0$$
{do not have common roots for generic $z\in \Bbb C^N$. 
 This contradicts our assumption that $g(y,z) = 0$ has repeated roots for 
generic $z$.}
\par
{We conclude from above that $i_j>1$ for some $j$. Now $f_j(y,z)$ is not a 
unit in ${\cal O}_0$, we have $f_j(0,0)=0$, so $f_j(0,z)\in (z_1, \cdots, 
z_N)$, that is, }
$$g(0,z) = g_m(z) = f_1^{i_1} \cdots f_h^{i_h} (0,z) \in (z_1, \cdots, z_N)^2$$
{because $i_j>1$. So the lemma is proved.}
\par
\vskip .2in
 \centerline{ \bf 2. The computations    }
\par
\vskip .1in
{We now start our computations. For simplicity of notations, we will assume 
that $q=(q_1,\cdots,q_s)=q_1=1$ in the rest of the paper.} 
\par
{Let $V\subset {\bf P^{n+r}}$ be a smooth complete intersection
 of degree $(d_1, d_2, \cdots, d_r)$, and} 
$$\beta_a = H^{n+1-a}\in H_{2(a-1)}(V, \Bbb C), \quad  
 1\leq a \leq n+1$$
{the $(n+1-a)$-th power of the 
hyperplane $H$, then $\{ \beta_a \}$ is a basis for the homology group 
$H_{inv}(V, \Bbb C)$. Any cohomology class $w\in H^*_{inv}(V, \Bbb C)$ 
can be written as }
$$w = t_1 \beta_1^* + t_2 \beta_2^* + \cdots + t_{n+1} \beta_{n+1}^*,$$
{here $\beta_a^* = H^{n+1-a}\in H^{2(n+1-a)}(V, \Bbb C)
 (1\leq a \leq n+1)$ in the ordinary cohomology. 
Then}
$$\Phi_{inv}^V(t_1,\cdots,t_{n+1})={d\over 6}\sum_{a+b+c=2n+3}t_at_b 
t_c+\sum_{k\geq 1}\sum_{\{k_a\}\in S_{V,k}} 
{{\sigma_k^V(k_1,\cdots,k_{n-1})t_1^{k_1}\cdots t_{n-1}^{k_{n-1}}} \over 
{k_1! \cdots k_{n-1}!}}e^{kt_n},$$
$$S_{V,k} = \{ \{k_a\}_{1\leq a\leq n-1} | \sum_{i=1}^{n-1} ik_{n-i} = 
(n+r+1-\sum_{i=1}^r d_i)k + n-3 > 0\},$$
$$(\beta_a^* \bullet_w \beta_b^*) (\beta_c) = {{\partial^3 \Phi^V}\over 
{\partial t_a \partial t_b \partial t_c}} (w),$$
{where $d=d_1d_2\cdots d_r$ is the degree of $V$, $\Phi^V$ is the 
Gromov-Witten prepotential, and 
$\sigma_k^V(k_1,\cdots,k_{n-1})$ is the number of degree $k$ rational 
curves in $V$ through $k_1$ points, $\cdots, k_{n-1}$ subspaces of 
dimension $n-2$ in general position.} 
\par
{Proposition 2, 3 in [B] and its higher degree analogy imply that there is 
a small positive number $\varepsilon = \varepsilon (d_1,\cdots,d_r,n) > 0$, 
such that $\Phi_{inv}^V(t_1,\cdots,t_{n+1})$ is well-defined and holomorphic 
in $t_1,\cdots,t_{n+1}$ when $|t_1|<\varepsilon, \cdots, |t_{n+1}| 
<\varepsilon $(also cf. [G],[J], [T]). 
\par 
{Since $X(w)\bullet_w$ is a linear operator on $H^*_{inv}(V, \Bbb C)$, its 
matrix $A(t_1,\cdots,t_{n+1})$ with respect to 
the basis $H^0=1, H^1, \cdots, H^n$ of $H^*_{inv}(V,\Bbb C)$ in the 
ordinary cohomology is a $(n+1)\times (n+1)$ matrix. Denote}
$$\eqalign{G(\lambda, t_1,\cdots,t_{n+1}) &= \det (\lambda I - A(t_1,\cdots,t_{n+1})),
\cr e &= 1+ \sum_{i=1}^r (d_i-1). \cr} $$
{In [B],[J],[G],[CJ], it was proved that }
$$G(\lambda, 0, \cdots, 0) = \lambda^{n+1} - d_1^{d_1}\cdots d_r^{d_r} 
\lambda^{e-1},$$
{when $e < n+1$ and $n\geq 3$. Therefore, we have}
\par
\vskip .1in
\noindent
{LEMMA 2. } {\it The only repeated 
root of the polynomial $G(\lambda, 0, \cdots, 0) = 0$ is $\lambda = 0$ when 
$e<n+1$ and $n\geq 3$.}
\par
\vskip .1in
\noindent
 
{Moreover, it is easy to see that $G(\lambda, 0, \cdots, 0)$ has distinct roots 
when $e=1, 2$ (cf. [KM],[T]), that is, $X(w)\bullet_{w}$ is semi-simple when $
e=1,2$ with $w=0$. Hence we  
may assume that $e\geq 3$.}
\par
{We choose}
$$w = t_{2e-3} \beta_{2e-3}^* = t_{2e-3} H^{n-2e+4} \in H_{inv}^*(V, \Bbb C)$$
{when $n>2e -3$, and assume that $t_{2e-3}$ is very close to 0. The rest of 
the computation will be made mod $(t_{2e-3})^2$. }
\par
{Consider}
$$(n-j)k_j + (n-m)k_m + (n-2e+3)k_{2e-3} = (n+2-e)k + (n-3).$$
{If $k_{2e-3} = 0$, then we have $k=1$ and $j+m=e+1$ because of the assumption 
$n>2e-3$. 
If $k_{2e-3}=1$, then we have either $k=1$ and $j+m=n-e+4$, or $k=2$ and 
$j+m=2$.}
\par
{Hence the $(n+1)\times (n+1)$-matrix for $H^1 \bullet_w$ with respect to 
the basis $H^0=1, H^1, \cdots, H^n$ of $H^*_{inv}(V,\Bbb C)$ in the 
ordinary cohomology is }
$$\pmatrix{0&1&0&\ldots&0&\ldots&0&\ldots&0\cr
0&0&1&\ldots&0&\ldots&0&\ldots&0\cr
\vdots&\vdots&\vdots&\ddots&\vdots&\ddots&\vdots&\ddots&0\cr
a_1t_{2e-3}&0&0&\ldots&1&\ldots&0&\ldots&0\cr
0&a_2t_{2e-3}&0&\ldots&0&\ldots&0&\ldots&0\cr
\vdots&\vdots&\vdots& &\vdots& &\vdots&\ddots&0\cr
b_1&0&0&\ldots&0&\ldots&1&\ldots&0\cr
0&b_2&0&\ldots&0&\ldots&0&\ldots&0\cr
\vdots&\vdots&\vdots& &\vdots&\ &\vdots&\ddots&1\cr
2c_1t_{2e-3}&0&0&\ldots&b_{e}&\ldots&a_{n-e+3}t_{2e-3}&\ldots&0\cr}.$$
{Here $e<n-e+3\leq n$ because of the assumption $n>2e-3$ and $e\geq 3$. When 
$e=3$, we have $a_1=a_{n-e+3}=a_n=0$ as $k_j$ is only defined for $1\leq j\leq n-1$.}
\par
{We next compute the matrix for $H^{n-2e+4}\bullet_w$.  
We only need to do the computation mod $(t_{2e-3})$ here.}
\par
{Consider}
$$(n-2e+3)k_{2e-3} + (n-j)k_j + (n-m)k_m = (n+2-e)k + (n-3).$$
{We choose $k_{2e-3}=1$ here. We have either $k=1$ and $j+m=n-e+4$, 
or $k=2$ and $j+m=2$. Hence the matrix for $H^{n-2e+4}\bullet_w$ is}
$$\pmatrix{0&0&\ldots&1&0&\ldots&0&\ldots&0\cr
0&0&\ldots&0&1&\ldots&0&\ldots&0\cr
\vdots&\vdots&\ddots&\vdots&\vdots&\ddots&\vdots&\ddots&\vdots\cr
a_1&0&\ldots&0&0&\ldots&1&\ldots&0\cr
0&a_2&\ldots&0&0&\ldots&0&\ldots&0\cr
\vdots&\vdots& &\vdots&\vdots& &\vdots&\ddots&\vdots\cr
c_1&0&\ldots&0&0&\ldots&a_{n-e+3}&\ldots&0\cr}.$$
{Here the 1 in the first row appears in the $(n-2e+5)$-th column.}
\par
{Now}
$$\eqalign{X(w) &= (n-2e+3)t_{2e-3} H^{n-2e+4} - (n-e+2) H^1,\cr
-{1\over {n-e+2}} X(w) &= H^1 - \delta t_{2e-3} H^{n-2e+4},\cr}$$
{here $\delta = {{n-2e+3}\over {n-e+2}}$, and $0<\delta <1$ because 
$n>2e-3$ and $e\geq 3$.}
\par
{Therefore the matrix for $(H^1 - \delta t_{2e-3} H^{n-2e+4})\bullet_w$ 
with respect to the basis $H^0=1, H^1, \cdots, H^n$ of $H^*_{inv}(V,\Bbb C)$ in the
 ordinary cohomology is }
$$A(t)=\pmatrix{0&1&0&\ldots&-t&0&\ldots&0&\ldots&0&\ldots&0\cr
0&0&1&\ldots&0&-t&\ldots&0&\ldots&0&\ldots&0\cr
\vdots&\vdots&\vdots& &\vdots&\vdots& &\vdots& &\vdots&\ddots&\vdots\cr
a_1^*t&0&0&\ldots&0&0&\ldots&1&\ldots&-t&\ldots&0\cr 
0&a_2^*t&0&\ldots&0&0&\ldots&0&\ldots&0&\ldots&0\cr
\vdots&\vdots&\vdots& &\vdots&\vdots&\ &\vdots& &\vdots&\ddots&\vdots\cr
b_1&0&0&\ldots&0&0&\ldots&0&\ldots&1&\ldots&0\cr
0&b_2&0&\ldots&0&0&\ldots&0&\ldots&0&\ldots&0\cr
\vdots&\vdots&\vdots& &\vdots&\vdots& &\vdots& &\vdots&\ddots&1\cr
c_1^*t&0&0&\ldots&0&0&\ldots&b_e&\ldots&a_{n-e+3}^*t&\ldots&0\cr},$$
{here $t=\delta t_{2e-3}$, $a_i^* = (\delta^{-1} - 1)a_i$, and $c_1^* = (2
\delta^{-1} - 1)c_1$.}
\par
{Now let $A_i(t)$ be the matrix obtained from $A(t)$ by replacing the 
$i-$th column of $A(t)$ by its derivative with respect to $t$.  
A detailed computation shows that}
$$\eqalign{\det A_1(0) &= (-1)^n (c_1^* - a_1^* b_e), \cr
\det A_{n-e+3}(0) &= (-1)^{n+1} b_1 (a_{n-e+3}^* + b_e), \cr
\det A_j(0) &= 0  \quad \quad \quad \hbox{ for } j\not= 1, n-e+3.\cr }$$
{As a result, }
$$\eqalign{ \hbox{det }A(t) &= (\det A_1(0) + \det A_2(0) + \cdots + \det A_{n+1}(0))
t\cr  
&= (-1)^n (c_1^* - a_1^*b_e - a_{n-e+3}^*b_1 - b_1b_e) t  
\quad \quad \quad \hbox{ mod }(t^2).\cr }$$
\par
{Following Beauville [B] we use $d l_j$ to denote the number of lines in 
$V$ meeting two general linear space of codimension $n-j$ and $n+1-e+j$ 
respectively.  
 We know that the varieties of lines and conics (with respect to 
the complex structure of $V$) contained 
in $V$ have the expected dimensions when $V$ is general [B]. }
\par
{In the case $e>3$, we have $b_1=b_e= l_0$. By proposition 2 in [B], 
we have}
$$a_1=a_{n-e+3} =  l_0.$$
{By the corollary of proposition 3 in [B], we have 
$c_1={1\over 2}l_0^2.$ Therefore }
$$\eqalign{c_1^* - a_1^*b_e - a_{n-e+3}^*b_1 - b_1b_e &= c_1^* - 2a_1^*b_1 - b_1^2\cr
   &=  -(\delta^{-1} - {1\over 2})l_0^2  \cr
   &< 0,\cr}$$
{because   
$0<\delta <1$. As a result, det $A(t) \not= 0$ mod $(t^2)$.}
\par
{In case $e=3$, $a_1=a_{n-e+3}=a_n=0$, we have }
$$\eqalign{c_1^* - a_1^*b_e - a_{n-e+3}^*b_1 - b_1b_e &= (\delta^{-1} - {3\over 2}) 
l_0^2 \cr
&= ({{n-1}\over {n-3}} - {3\over 2}) l_0^2 \cr &\not= 0\cr}$$
{except $n=7$ and $e=3$. Again, we have $\det A(t) \not= 0$ mod $(t^2)$.}
\par
{Finally, we give}
\par
\vskip .1in
\noindent
{\it Proof of Theorem 1.} { By the above computation, we have }
$$\eqalign{G(0,\cdots,0,t_{2e-3},0,\cdots,0) &= \det ( - 
A(0,\cdots,0,t_{2e-3},0,\cdots,0)) \cr 
&= \det ( - A(t)) \cr
&\not= 0 \quad \quad \hbox{mod }(t^2_{2e-3}).\cr }$$
{Therefore, 
Lemma 1 and 2 implies that the polynomial 
$G(\lambda,0,\cdots,0,t_{2e-3},0,\cdots,0) = 0$ has distinct roots for 
generic $t_{2e-3}\in \Bbb C$. However, being semi-simple is an open 
conditions, we conclude that the polynomial } 
$$G(\lambda, t_1,\cdots,t_{n+1}) = 0$$
{has distinct roots for generic 
$(t_1,\cdots,t_{n+1})\in \Bbb C^{n+1}$, that is, $X(w)\bullet_w$ is 
semi-simple for generic $w$.}
\par
\vskip .1in
 
\par

\par
\centerline{\bf References}
\par
\vskip .1in
\par
\noindent
\item{[B]} A. Beauville, {\it Quantum cohomology of complete intersections}, 
 alg-geom/9501008. 
\par
\noindent
\item{[CJ]} A. Collino and M. Jinzenji, {\it On the structure of small quantum 
cohomology rings for projective hypersurfaces,}, preprint, 1996.
\par
 \noindent
\item{[D]} B. Dubrovin, {\it Integrable systems in topological field theory}, 
Nucl. Phys. {\bf B379} (1992), 627-689. Also {\it Geometry of 2d topological 
field theory}, preprint 1994. 
\par
\noindent
\item{[G]} A. Givental, {\it Equivariant Gromov-Witten invariants,} preprint, 1996.
\par 
\noindent
\item{[GH]} P. Griffiths and J. Harris, 
 {\it Principles of Algebraic Geometry}, Wiley, New York, 1978. 
\par
\noindent
\item{[J]} M. Jinzenji, {\it On quantum cohomology rings for hypersurfaces in 
${\Bbb C}P_{N-1}$}, preprint, hep-th/9511206.  
\par
 \noindent
\item{[KM]} M. Kontsevich and Yu. Manin, {\it Gromov-Witten classes, 
quantum cohomology, and enumerative geometry}, Comm. Math. Phy. {\bf 164} 
(1994), 525-562. 
\par
\noindent
\item{[MS]} D. McDuff and D. Salamon, {\it J-holomorphic curves and 
quantum cohomology}, AMS, University Lectures Series, v. 6, 1994. 
\par
\noindent
\item{[RT]} Y. Ruan and G. Tian, {\it A mathematical theory of quantum 
cohomology}, J. Diff. Geom. {\bf 42} (1995), 259-367. 
\par
\noindent
\item{[T]} G. Tian, {\it Quantum cohomology and its associativity}, Proc. 
of 1st Current developments in Mathematics, Cambridge, 1995.  
\par
\vskip .1in
\par
Gang Tian
\par
Department of Mathematics
\par
Massachusetts Institute of Technology
\par
Cambridge, MA 02139
\par
e-mail: tian@math.mit.edu
\par
\vskip .1in
Geng Xu
\par
Department of Mathematics
\par
Johns Hopkins University
\par
Baltimore, MD 21218
\par
e-mail: geng@math.jhu.edu 

\end